\documentclass{amsart}
\usepackage{CJKutf8}
\usepackage{amsmath,amsfonts,amssymb,amscd,bm,latexsym}
\usepackage[colorlinks=true]{hyperref}
\usepackage{tikz}
\usepackage{pgflibraryarrows}
\usepackage{pgflibrarysnakes}
\usepackage{pgflibraryarrows}
\usepackage{pgflibrarysnakes}
\hypersetup{urlcolor=blue, citecolor=blue}
\numberwithin{equation}{section} 

\def\<{\langle}             \def\>{\rangle}

\newtheorem{thm}{Theorem}[section]

\newtheorem{lem}[thm]{Lemma}

\theoremstyle{definition}

\newtheorem{rem}{Remark}[section]

\newcommand{\beeq}{\begin{equation}}\newcommand{\eneq}{\end{equation}}
    \newcommand{\be}{\beta}
\newcommand{\de}{\delta}    
\newcommand{\vep}{\epsilon}\newcommand{\ep}{\varepsilon}
  
    \newcommand{\la}{\lambda}

\newcommand{\ga}{\gamma}    \newcommand{\Ga}{\Gamma}
\newcommand{\R}{\mathbb{R}}

\newenvironment{prf}{\noindent {\bf Proof.} }{\endprf\par}
\def \endprf{\hfill  {\vrule height6pt width6pt depth0pt}\medskip}
\newcommand{\pa}{\partial}
\newcommand{\les}{{\lesssim}}
\newcommand{\gt}{{\gtrsim}}
\newcommand{\supp}{\,\mathop{\!\mathrm{supp}}}

\newcommand{\CO}{\mathcal{O}}

\numberwithin{equation}{section}
\title[Wave system with damping and potential]
      {Blow up and lifespan estimates for systems of semi-linear wave equations with damping and potential }
\author{Mengliang Liu}
\address{Department of Mathematics\\                	
Guangxi University\\                Nanning 530000, P. R. China}
\email{liumengliang@alu.gxu.edu.cn}

\date{\today}
\dedicatory{} \commby{}

\begin{document}

\begin{abstract}
In this paper, we consider the semi-linear wave systems with power-nonlinearities and a large class of space-dependent damping and potential. We obtain the same blow-up regions and the lifespan estimates for three types wave systems, compared with the systems without damping and potential.
\end{abstract}



\maketitle

\section{Introduction}
In this paper, we study the finite time blow-up phenomenon of three kinds of semi-linear wave systems with space dependent damping and potential. More precisely, we consider the following small data Cauchy problem with power-nonlinearities
\beeq
\label{1}
\begin{cases}
u_{tt}-\Delta u+D_1(x)u_t+V_{1}(x)u=N_{1}(v,v_t), &(t,x)\in (0,T) \times \R^{n}, \\
v_{tt}-\Delta v+D_2(x)v_t+V_{2}(x)v=N_{2}(u,u_t), &(t,x)\in  (0,T) \times \R^{n},\\
u(0, x)=\ep u_0(x), \ \ \ u_t(0, x)=\ep u_1(x),\ \ \ &x \in \R^{n},\\
v(0, x)=\ep v_0(x), \ \ \ v_t(0, x)=\ep v_1(x),\ \ \ &x \in \R^{n}\ .\\
\end{cases}
\eneq
Here $(u_0, u_1), (v_0, v_1)\in C^\infty_c(\R^n)$,  and the small parameter $\ep>0$ measures the size of the data. As usual, to show blow up, we assume  both $(u_0, u_1)$ and $(v_0, v_1)$ are nontrivial, nonnegative and supported in $B_R:= \{x\in\R^n: r\le R\}$ for some $R>0$, where $|x|=r$.  For the coefficient of dampings $D_i(x)$ and potentials $V_{i}(x)$, we assume $V_{i}(x)$, $D_{i}(x) \in C(\R^{n}) \cap C^{\delta}(B_{\delta})$ for some $\de>0$. In addition, we assume 
\begin{align}
\label{4}
0\leq D_{i}(x)=D_{i}(|x|)\leq \alpha_{1}(1+|x|)^{-\be}, \ 0\leq \alpha_{1}\in \R, 1 <\be\in \R, i=1, 2.
\end{align}
\begin{align}
\label{3}
0\leq V_{i}(x)=V_{i}(|x|)\leq \alpha_{2}(1+|x|)^{-\omega}, \ 0\leq\alpha_{2}\in \R, 2<\omega\in \R,i=1,2.
\end{align}
In order to motivate the study of \eqref{1}, let us recall some semi-linear models which are strongly related to this weakly coupled system. We begin with the following system
\begin{align}
\label{2}
\begin{cases}
u_{tt}-\Delta u=G_{1}(v,v_{t}), \ \ \ \ \   t>0, x\in\R^{n}, \\
v_{tt}-\Delta v=G_{2}(u,u_{t}), \ \ \ \ \   t>0, x\in\R^{n}, \\
u(0, x)=\ep u_0(x), \ u_t(0, x)=\ep u_1(x)\\
v(0, x)=\ep v_0(x), \ v_t(0, x)=\ep v_1(x)\\
\end{cases}
\end{align}
\ \ \ When $ G_{1}(v,v_{t})=|v|^{p}$, $G_{2}(u,u_{t})=|u|^{q}$, Santo-Georgiev-Mitidieri \cite{S-D-G-V-E} proved there exist a critical curve in $ (p,q) $-plane, 
$$\Ga_{SS}(n,p,q)=\max\{\frac{p+2+q^{-1}}{pq-1}, \frac{q+2+p^{-1}}{pq-1}\}-\frac{n-1}{2} .$$
They showed that when $ \Ga_{SS}(n,p,q)<0$, $ n\geq2$, the system \eqref{2} has a global in time solution for sufficiently small $\ep$, while a solution for some positive initial data blows up in finite time if $\Ga_{SS}(n,p,q)>0$. Palmieri-Takamura \cite{P-H} obtained the upper bound of lifespan estimate:
\begin{align}
\label{thm1up}
 T_\ep\lesssim
 \begin{cases}
 \ep^{-\Ga_{SS}^{-1}(n,p,q)}, \ \ &\Ga_{SS}(n,p,q)>0,\\
 \exp\big(\ep^{-\min\{p(pq-1),q(pq-1)\} }\big),\  \ &\Ga_{SS}(n,p,q)=0,p\neq q,\\
 \exp\big(\ep^{-p(p-1)}\big),\  \ &\Ga_{SS}(n,p,q)=0, p=q=p_{S}(n).
 \end{cases}
\end{align}
See also \cite{S-E}  \cite{H-K-O} \cite{I-M-S-M-W-K} for relevant works.

When $G_{1}(v,v_{t})=|v_{t}|^{p}$, $G_{2}(u,u_{t})=|u_{t}|^{q}$, it has been investigated by Deng \cite{D-K}, Ikeda-Sobojima-Wakasa \cite{I-M-S-M-W-K}, Xu \cite{X-W}. Where, Deng \cite{D-K} has obtained the blow up region of system \eqref{2} in 
\begin{align*}
1<pq<\infty
\begin{cases}
1\leq p<\infty, \ \ \ &n=1,\\
1<p\leq2,\ \ \ &n=2,\\
p=1,\ \ \ &n=3.
\end{cases}
\end{align*}
\begin{align*}
1<pq\leq\frac{(n+1)p}{(n-1)p-2}
\begin{cases}
2<p\leq3, \ \ \ &n=2,\\
1<p\leq2,\ \ \ &n=3,\\
1<p\leq\frac{n+1}{n-1},&n\geq4.
\end{cases}
\end{align*}
Ikeda-Sobajima-Wakasa \cite{I-M-S-M-W-K} obtained the finite time blow up results when
$$\Ga_{GG}(p,q,n):=\max\{\frac{p+1}{pq-1}, \frac{q+1}{pq-1}\}-\frac{n-1}{2} \geq0,\   p, q>1,\ n \geq 2\ ,$$
and obtain the upper bound estimates:
\begin{align}
\label{thmup2}
T_{\ep}\lesssim
\begin{cases}
\ep^{-\Ga^{-1}_{GG}(n,p,q)},\ \ &\Ga_{GG}(n,p,q)>0,\\
\exp(\ep^{-(pq-1)}),\ \ &\Ga_{GG}(n,p,q)=0,\ \ \ p\neq q,\\
\exp(\ep^{-(p-1)}),\ \ &\Ga_{GG}(n,p,q)=0,\ \ \ p=q.\\
\end{cases}
\end{align}

When $G_{1}(v,v_{t})=|v|^{q}$, $G_{2}(u,u_{t})=|u_{t}|^{p}$, Hidano-Yokoyama \cite{H-K-Y-K} showed the upper bound of lifespan estimates:
\begin{align}
T_{\ep}\lesssim\ep^{-\frac{p(pq-1)}{p+2-(\frac{(n-1)p}{2}-1)(pq-1)}},
\end{align}
when $n\geq2$, $1<q$, $1<p<\frac{2n}{n-1}$ and $(\frac{(n-1)p}{2}-1)(pq-1)<p+2$.
Ikeda-Sobojima-Wakasa \cite{I-M-S-M-W-K} introduced two kind of exponent for the system (\ref{2}),
$$F_{SG,1}(n,p,q)=(\frac{1}{p}+1+q)(pq-1)^{-1}-\frac{n-1}{2},$$
$$F_{SG,2}(n,p,q)=(\frac{1}{q}+2)(pq-1)^{-1}-\frac{n-1}{2}.$$
They showed the solution for some positive data blows up in finite time if 
$$\Ga_{SG}(n,p,q)=\max\{F_{SG,1}(n,p,q),F_{SG,2}(n,p,q)\}\geq0\ ,$$
and obtained the upper bound estimates:
\begin{align}
\label{thmup3}
T_{\ep}\lesssim
\begin{cases}
\ep^{-\Ga^{-1}_{SG}(n,p,q)},\ \ &\Ga_{SG}(n,p,q)>0,\\
\exp(\ep^{-q(pq-1)}),\ \ &F_{SG,1}(n,p,q)=0>F_{SG,2}(n,p,q),\\
\exp(\ep^{-q(pq-1)}),\ \ &F_{SG,2}(n,p,q)=0>F_{SG,1}(n,p,q),\\
\exp(\ep^{-(pq-1)}),\ \ &F_{SG,2}(n,p,q)=0=F_{SG,1}(n,p,q).
\end{cases}
\end{align}

Concerning the influence of dampings and potentials to the blow up region of semilinear wave equations, Lai-Tu \cite{L-T} considered the small data Cauchy problems
\begin{align}
\begin{cases}
u_{tt}-\Delta u+\frac{\mu}{(1+|x|)^{\be}}u_{t}=N(u, u_t), (t,x) \in(0, T)\times \R^n, \\
u(0, x)=\ep f(x), \ u_t(0, x)=\ep g(x) \\
\end{cases}
\end{align}
where $\mu>0,$ $\be>2$ are constants. They showed that such damping terms does not effect the upper bound of lifespan estimates:
\begin{align*}
N(u, u_t)=|u|^p\ , \ T_{\ep}\lesssim
\begin{cases}
\ep^{-\frac{2(p-1)}{n+1-(n-1)p}},\ \ 1<p\leq\frac{n}{n-1},\\
\ep^{-\frac{2p(p-1)}{\ga(p,n)}},\ \ \frac{n}{n-1}<p<p_{S}(n).
\end{cases}
\end{align*}
\begin{align*}
N(u, u_t)=|u_t|^p\ , \ 
T_{\ep}\lesssim
\begin{cases}
\ep^{-(\frac{1}{p-1}-\frac{n-1}{2})^{-1}},\ \ 1<p\leq p_{G}(n),\\
\exp\big(\ep^{-(p-1)}\big),\ \ \ \ \ \ \ \ \  p=p_{G}(n).
\end{cases}
\end{align*}
Here, $p_{S}(n)$ is the Strauss exponent (see \cite{W-A}) and $p_S$ is the positive root of the quadratic equation
$$
2+(n+1)p-(n-1)p^{2}=0.
$$
And $p_G(n)=1+2/(n-1)$ is Glassey exponent (see \cite{R-T}). Recently, 
Lai-Liu-Tu-Wang \cite{L-L-W-T} considered the following semi-linear wave equation
\begin{align}
\label{lmy}
\begin{cases}
u_{tt}-\Delta u+D(x)u_{t}+V(x)u=|u|^{p}, \ \ \ \ \ (t, x)\in  (0, T)\times \R^n, \\
u(0, x)=\ep f(x), \ u_t(0, x)=\ep g(x)
\end{cases}
\end{align}
where $n\geq2$. They showed the problem exists a large class of damping and potential functions of critical/long range in the blow up part. In particular, they showed that for the coefficient of dampings and potentials satisfying \eqref{4}, \eqref{3}, the solution of \eqref{lmy} will blow up under the Strauss exponent $p_S(n)$. Based on this work, we can imagine such damping and potential will not effect the blow up region of system \eqref{1}.  

In this work, we obtain the blow up results and the upper bound lifespan estimates to three types of system \eqref{1} under the exponent $\Ga_{SS}$, $\Ga_{GG}$, $\Ga_{SG}$ with a large class of damping and potential functions \eqref{4} \eqref{3}.

Before proceeding, we give a definition of the weak solution.

{\bf Definition of weak solutions:} Suppose that $u, v\in C(([0,T), H^{1}(\R^n))\cap C^{1}([0,T),L^{2}(\R^n))$ and $u\in\ L^{q}_{loc}((0,T) \times \R^n),$ $v\in\ L^{p}_{loc}((0,T) \times \R^n)$ for Theorem \ref{liu1}, $u_{t}\in\ L^{q}_{loc}((0,T)\times \R^n),$ $v_{t}\in\ L^{p}_{loc}((0,T) \times \R^n)$ for Theorem \ref{liu2}, $u_{t}\in\ L^{p}_{loc}((0,T) \times \R^n),$ $v\in\ L^{q}_{loc}((0,T) \times \R^n)$ for Theorem \ref{liu3} and 
\begin{align} 
\nonumber
&\ep\int_{\R^{n}}u_{1}(x)\Psi(0,x)dx+\ep\int_{\R^{n}}u_{0}(x)D_{1}(x)\Psi(0,x)dx
+\int_{0}^{T}\int_{\R^{n}}N_{1}(v,v_{t})\Psi(t,x)dxdt\\ \nonumber
&\ \ \ \ \ \ \ \ \ \ \ \ \  =-\int_{0}^{T}\int_{\R^{n}}u_{t}(t,x)\Psi_{t}(t,x)dxdt+\int_{0}^{T}\int_{\R^{n}}\nabla u(t,x)\nabla \Psi(t,x)dxdt\\ \label{l1}
&\ \ \ \ \ \ \ \ \ \ \ \ \ -\int_{0}^{T}\int_{\R^{n}}u(t,x)D_{1}(x)\Psi_{t}(t,x)dxdt+\int_{0}^{T}\int_{\R^{n}}V_{1}(x)u(t,x)\Psi(t,x)dxdt
\end{align}
\begin{align} \nonumber
\label{l2}
&\ep\int_{\R^{n}}v_{1}(x)\Psi(0,x)dx+\ep\int_{\R^{n}}v_{0}(x)D_{2}(x)\Psi(0,x)dx
+\int_{0}^{T}\int_{\R^{n}}N_{2}(u,u_{t})\Psi(t,x)dxdt\\ \nonumber
&\ \ \ \ \ \ \ \ \ \ \ \ \ =-\int_{0}^{T}\int_{\R^{n}}v_{t}(t,x)\Psi_{t}(t,x)dxdt+\int_{0}^{T}\int_{\R^{n}}\nabla v(t,x)\nabla \Psi(t,x) dxdt\\
&\ \ \ \ \ \ \ \ \ \ \ \ \ \ -\int_{0}^{T}\int_{\R^{n}}v(t,x)D_{2}(x)\Psi_{t}(t,x)dxdt+\int_{0}^{T}\int_{\R^{n}}V_{2}(x)v(t,x)\Psi(t,x)dxdt\\ \nonumber
\end{align}
for any $ \Psi(t,x)\in C_{0}^{\infty}([0,T)\times\R^n)$. The supremum of all such time of existence,  $T$, is called to be the
lifespan to the problem \eqref{1}, denoted by $T_\vep$.

As usual, to show blow up, we need the solution satisfy finite speed of propagation. We suppose that the weak solution satisfies
\begin{equation}
\supp\ (u,v) \subset  \{(x,t)\in \R^n\times\ [0,T) : |x|\leq t+R \}.
\end{equation}
The main results in this paper are illustrated by the following Theorems.
\begin{thm} 
\label{liu1}
Let $n\geq3$. Consider the semi-linear wave systems problem \eqref{1} with nonlinear terms
$N_{1}(v,v_t)=|v|^{p}$, $N_{2}(u,u_t)=|u|^{q}$, damping \eqref{4}, potential \eqref{3}. Any weak solution with nontrivial, nonnegative, compactly supported data will blow up in finite time. In addition, there exists a constant $\ep_0$ such that for any $\ep\in (0, \ep_0)$, the lifespan $T_\ep$ satisfies \eqref{thm1up}. 
 \end{thm}
\begin{thm} 
\label{liu2}
Let $n\geq3$. Consider the semi-linear wave systems problem \eqref{1} with nonlinear terms
$N_{1}(v,v_t)=|v_{t}|^{p}$, $N_{2}(u,u_t)=|u_{t}|^{q}$, damping \eqref{4}, potential \eqref{3}. Any weak solution with nontrivial, nonnegative, compactly supported data will blow up in finite time. In addition, there exists a constant $\ep_0$ such that for any $\ep\in (0, \ep_0)$, the lifespan $T_\ep$ satisfies \eqref{thmup2}. 
\end{thm}

\begin{thm}
\label{liu3} 
Let $n\geq3$. Consider the semi-linear wave systems problem \eqref{1} with nonlinear terms $N_{1}(v,v_t)=|v|^{q}$, $N_{2}(u,u_t)=|u_{t}|^{p}$, damping \eqref{4} and potential \eqref{3}.  Any weak solution with nontrivial, nonnegative, compactly supported data will blow up in finite time. In addition, there exists a constant $\ep_0$ such that for any $\ep\in (0, \ep_0)$, the lifespan $T_\ep$ satisfies \eqref{thmup3}. 
\end{thm}
\begin{rem}
When damping and potential fuctions satisfying \eqref{4}, \eqref{3}, we conjecture the critical curve of the solution to system \eqref{1} is the same as  system \eqref {2}(with corresponding nonlinear term). That is to say, such damping and potential functions does not effect the critical curve of system \eqref{1}.
\end{rem}

\section{PRELIMINARY}

In this section, we collect some Lemmas we will use later. For the strategy of the proof, we use test function methods, see, e.g., \cite{Q-Z}, \cite{B-Q-Z}, \cite{I-M-S-M-W-K}.

\begin{lem}[{\bfseries Lemma 2.2 in \cite{L-L-W-T}}]
\label{le5}
Support $0\leq V \in C(\R^{n})$, $V=V=\CO(\<r\>^{-\omega})$ for some $\omega>2$, and $V \in C^{\delta}(B_{\delta})$ for some $\delta>0$. In addition, assume $V$ is nontrivial when $n=2$. Then there exist a $C^{2}$ solution of
\begin{equation}
\label{lp3}
\Delta \phi_0=V\phi_0,\ \  x \in R^{n}, n\geq2,
\end{equation} 
satisfying
\begin{align}
\label{lp4}
\phi_0\simeq
\begin{cases}
\ln(r+2), &n=2,\\
1, &n\geq3.
\end{cases}
\end{align}
Moreover, when $n=2$ and $V=0$, it is clear that $\phi_0=1$ is a solution \eqref{lp3}.
\end{lem}

\begin{lem}[{\bfseries Lemma 2.4 in \cite{L-L-W-T}}]
\label{lp5}
Let $V=V(r)$, $D=D(r) \in C(\R^{n})\cap C^{\delta}(B_{\delta})$ for some $\delta>0$. Suppose $V\geq0$, $D\geq-\la-\frac{V}{\la}$, and also that for some $d_{\infty} \in \R$, $R>1$ and $r\geq R$, we have 
\begin{equation}
\label{lp6}
V(r) \in L^{1}([R, \infty)),\ \ D(r)=\frac{d_{\infty}}{r}+D_{\infty}(r),\ \ D_{\infty} \in L^{1}([R, \infty)).
\end{equation}
Then there exists a $C^{2}$ solution of 
\begin{equation}
\label{lp7}
\Delta \phi_{\la}=(\la^{2}+\la D+V)\phi_{\la},\ \  \la>0,\ \ x \in \R^{n},
\end{equation}
satisfying 
\begin{equation}
\phi_{\la}\simeq\<r\>^{-\frac{n-1-d_{\infty}}{2}}e^{\la r}.
\end{equation}
\end{lem}

\begin{lem}
\label{le3}
Concerning \eqref{4} \eqref{3}, we have $d_\infty=0$ in Lemma \ref{lp5}. 
Let $a>0$. Define
\begin{equation*}
b_{a}(t,x)=\int^{1}_{0}e^{-\la t}\phi_{\la}(x)\la^{a-1}d\la, \ x\in \R^{n}
\end{equation*}
$ b_{a}(t,x) $ lie on the facts that they satisfy
\begin{equation}
\label{lll1}
\pa^{2}_{t}b_{a}-\Delta b_{a}-D_{i}(x)\pa_{t}b_{a}+V_i(x)b_a=0, \ \ i=1 , 2\ .
\end{equation}
\begin{equation*}
\pa_{t}b_{a}=-b_{a+1}\ .
\end{equation*}
and enjoy the asymptotic behavior for $n\geq 2$ and $ r\leq t+R $
\begin{equation}
\label{j1}
b_{a}(t,x)\gtrsim (t+R)^{-a}, a>0,
\end{equation}
\begin{align}
 b_{a}(t,x)\lesssim
 \begin{cases}
 \label{j2}
 \ (t+R)^{-a}, \     &\  0<a<\frac{n-1}{2},\\
\ (t+R)^{-\frac{n-1}{2}}(t+R+1-|x|)^{\frac{n-1}{2}-a}, \ &\ a>\frac{n-1}{2}.
\end{cases}
\end{align}
\end{lem}
\begin{prf}See the proof of (6.7) (6.8) in \cite{L-L-W-T}. 
\end{prf}
\begin{lem}[{\bfseries Lemma 3.1 in \cite{L-T}}]
\label{le1}
If $\be>0,$ then for any $\alpha \in \R$ and a fixed constant R, we have
\begin{align*}
\int^{t+R}_{0}(1+r)^{\alpha}e^{-\be(t-r)}dr \ \les \ (t+R)^{\alpha}\ .
\end{align*}
\end{lem}

Finally, we need the following ODE Lemma to show the finite time blow up in some critical cases. 
\begin{lem}[{\bfseries Lemma 3.10 in \cite{I-M-S-M-W-K}}]
\label{le4}
Let $2<t_{0}<T $, $0\leq\phi\in C^{1}([t_{0},T))$. Assume that
\begin{align*}
\begin{cases}
\delta\leq K_{1}t\phi'(t),\  t\in(t_{0},T)\\
\phi(t)^{p_{1}}\leq K_{2}t(\ln t)^{p_{2}-1}\phi'(t), \  t\in(t_{0},T)
\end{cases}
\end{align*}
with $\delta$, $K_{1}$, $K_{2}>0$ and $p_{1}$, $p_{2}>1$.  If $p_{2}<p_{1}+1$, then there exists positive constant $\delta_{0}$ and $K_{3}$ (independent of $\delta$) such that
\begin{align*}
T\leq\exp(K_{3}\delta^{-\frac{p_{1}-1}{p_{1}-p_{2}+1}})\nonumber
\end{align*} 
for all $\delta \in (0,\delta_{0})$.
\end{lem}
{\bf{Set up:}}  We take $\la=1$ in Lemma \ref{lp5} and $\phi(x)=\phi_1(x)$,    $\Phi(t,x)=e^{-t}\phi(x)$, it is easy to see that $\Phi$ is the solution of Linear wave equation
$$\pa^{2}_{t}\Phi-\Delta\Phi-D_{i}(x)\pa_{t}\Phi+V_i\Phi=0,   \ i=1, 2. $$

\section{Proof of Theorems \ref{liu1}-\ref{liu3}}\label{proof}
In this section, we give the proof of Theorems \ref{liu1}-\ref{liu3}.
\subsection{Proof of Theorem \ref{liu1}}
\subsubsection{{\bf Blow up region $\Ga_{SS}(n,p,q)>0$}}
Let $ \eta(t)\in C^{\infty} ([0,\infty))$ satisfies\\
\begin{align}
\eta(t)=
\begin{cases}
\ 1,       &0<t\leqslant\frac{1}{2},\\
\ decreasing,   &\frac{1}{2}<t<1,\\
\ 0,             &t\geqslant1,
\end{cases}
\end{align}
and
$$\eta_{M}(t)=\eta(\frac{t}{M}),\ M\in(1,T_{\ep}).$$
In addition, we take $T\in [M, T_\ep)$. 

{\bf {(\uppercase\expandafter{\romannumeral1}) Test function 1: $\Psi(t,x)=\eta^{2p^{\prime}}_{M}(t)\phi_0(x)$}}

Firstly, we choose the test function $\Psi(t,x)=\eta^{2p^{\prime}}_{M}(t)\phi_0(x)$ to replace $\Psi(t,x)$ in \eqref{l1} and applying Lemma \ref{le5}, we get that
\begin{align*} \nonumber
&\ep\int_{\R^{n}}u_{1}(x)\phi_0 dx+\ep\int_{\R^{n}}u_{0}(x)D_{1}(x)\phi_0 dx
+\int_{0}^{M}\int_{\R^{n}}|v|^{p}\eta_{M}^{2p^{\prime}}\phi_0 dxdt\\  \nonumber
&=\int_{0}^{M}\int_{\R^{n}}u\phi_0\big(\pa^{2}_t \eta_{M}^{2p'}-D_{1}(x)\pa_t \eta_{M}^{2p'}\big)dxdt-\int_{0}^{M}\int_{\R^{n}}u\eta_{M}^{2p'}(\Delta \phi_0-V_{1}(x)\phi_0)dxdt\\  
&=\int_{0}^{M}\int_{\R^{n}}u\phi_0\big(\pa^{2}_t \eta_{M}^{2p'}-D_{1}(x)\pa_t \eta_{M}^{2p'}\big)dxdt\\
&=I_{1}+I_{2}
\end{align*}
Noting that when $p\leq q$ we have $2p'q-2q\geq 2p'$.  We apply H\"older inequality and \eqref{lp4} to the following estimates
\begin{align} \nonumber
\label{l3}
|I_{1}|\lesssim&M^{-2}\int_{0}^{M}\int_{\R^{n}}u(t,x)\theta^{2p'-2}_{M}dxdt\\\nonumber
\lesssim&M^{-2}\Big(\int_{0}^{M}\int_{\R^{n}}|u|^{q}\theta^{(2p'-2)q}_{M}dxdt\Big)^{\frac{1}{q}}\Big(\int_{\frac{M}{2}}^{M}\int_{\R^{n}}dxdt\Big)^{\frac{1}{q'}}\\ 
\lesssim&M^{-2}\Big(\int_{0}^{M}\int_{\R^{n}}|u|^{q}\theta^{2p'}_{M}dxdt\Big)^{\frac{1}{q}}\Big(\int_{\frac{M}{2}}^{M}\int_{\R^{n}}dxdt\Big)^{\frac{1}{q'}} \   \ \ \ \ \ \ \ (p\leq q)\\ \nonumber
\lesssim&M^{-2}\Big(\int_{0}^{M}\int_{\R^{n}}|u|^{q}\theta^{2p'}_{M}dxdt\Big)^{\frac{1}{q}}\Big(\int_{\frac{M}{2}}^{M}\int_{0}^{t+R}r^{n-1}drdt\Big)^{\frac{1}{q'}}\\ \nonumber
\lesssim& M^{\frac{nq-n-q-1}{q}}\Big(\int^{M}_{0}\int_{\R^{n}}|u|^{q}\theta^{2p'}_{M}dxdt\Big)^{\frac{1}{q}},
\end{align}
\begin{align} \nonumber
\label{l4}
|I_{2}|\lesssim&M^{-1}\int_{0}^{M}\int_{\R^{n}}u(t,x)D_{1}(x)\theta^{2p'-1}_{M}dxdt\\  \nonumber
\lesssim&M^{-1}\Big(\int_{0}^{M}\int_{\R^{n}}|u|^{q}\theta^{(2p'-1)q}_{M}dxdt\Big)^{\frac{1}{q}}\Big(\int_{\frac{M}{2}}^{M}\int_{\R^{n}}(D_{1}(x))^{q'}dxdt\Big)^{\frac{1}{q'}}\\  
\lesssim&M^{-1}\Big(\int_{0}^{M}\int_{\R^{n}}|u|^{q}\theta^{2p'}_{M}dxdt\Big)^{\frac{1}{q}}\Big(\int_{\frac{M}{2}}^{M}\int_{0}^{t+R}(1+r)^{-\be q'}r^{n-1}drdt\Big)^{\frac{1}{q'}}\\ 
\nonumber
\lesssim&M^{-1}\Big(\int^{M}_{0}\int_{\R^{n}}|u|^{q}\theta^{2p'}_{M}dxdt\Big)^{\frac{1}{q}}
\begin{cases} \nonumber
\ (M\ln(M))^{\frac{1}{q'}}, \ \ \  &n-\be q'=0,\\
\ M^{(n+1-\be q')\frac{1}{q'}},\ \ \  &n-\be q'>0,\\
\ M^{\frac{1}{q'}},\ \ \  &n-\be q'<0.
\end{cases}\\ \nonumber
\lesssim&M^{\frac{nq-n-\be q-1}{q}}\Big(\int^{M}_{0}\int_{\R^{n}}|u|^{q}\theta^{2p'}_{M} dxdt\Big)^{\frac{1}{q}}\\ \nonumber
\lesssim& M^{\frac{nq-n-q-1}{q}}\Big(\int^{M}_{0}\int_{\R^{n}}|u|^{q}\theta^{2p'}_{M}dxdt\Big)^{\frac{1}{q}}.\nonumber
\end{align} 
By combining \eqref{l3} and \eqref{l4}, we get that
\begin{equation}
\ep C(u_{0},u_{1})+\int_{0}^{M}\int_{\R^{n}}|v|^{p}\eta_{M}^{2p^{\prime}}dxdt\lesssim M^{\frac{nq-n-q-1}{q}}\Big(\int^{M}_{0}\int_{\R^{n}}|u|^{q}\theta^{2p'}_{M}dxdt\Big)^{\frac{1}{q}},
\end{equation}
then we get that
\begin{align}
\label{m1}
\int_{0}^{M}\int_{\R^{n}}|v|^{p}\eta_{M}^{2p^{\prime}}dxdt\lesssim M^{\frac{nq-n-q-1}{q}}\Big(\int^{M}_{0}\int_{\R^{n}}|u|^{q}\theta^{2p'}_{M}dxdt\Big)^{\frac{1}{q}}.
\end{align}
By a similar way, we choose the test function with $\Psi(t,x)=\eta^{2p^{\prime}}_{M}(t)\phi_0(x)$ to replace $\Psi(t,x)$ in \eqref{l2}, we can get
\begin{equation}
\label{ll}
\int_{0}^{M}\int_{\R^{n}}|u|^{q}\eta_{M}^{2p^{\prime}}dxdt\lesssim M^{\frac{np-n-p-1}{p}}\Big(\int^{M}_{0}\int_{\R^{n}}|v|^{p}\theta^{2p'}_{M}dxdt\Big)^{\frac{1}{p}}.
\end{equation}
Applying \eqref{m1} and \eqref{ll} yields
\begin{equation}
\label{l6}
\int_{0}^{M}\int_{\R^{n}}|v|^{p}\eta_{M}^{2p^{\prime}}dxdt\lesssim M^{\frac{nqp-2p-n-pq-1}{qp-1}}, 
\end{equation}
\begin{equation}
\label{l7}
\int_{0}^{M}\int_{\R^{n}}|u|^{q}\eta_{M}^{2p^{\prime}}dxdt\lesssim M^{\frac{nqp-2q-n-pq-1}{qp-1}}.
\end{equation}

{\bf {(\uppercase\expandafter{\romannumeral2}) Test function 2: $\Psi(t,x)=\eta^{2p'}_{M}(t)e^{-t}\phi(x)=\eta^{2p'}_{M}(t)\Phi(t,x)$}}

Secondly, we choose another test function $\Psi=\eta^{2p'}_{M}\Phi(t,x)$ to replace $\Psi(t,x)$ in \eqref{l1} and applying Lemma \ref{lp5}, we get that
\begin{align*}
&\ep C(u_{0},u_{1})+\int_{0}^{M}\int_{\R^{n}}|v|^{p}\eta_{M}^{2p^{\prime}}\Phi dxdt\\
=&-\int_{0}^{M}\int_{\R^{n}}u_{t}\pa_{t} (\eta_{M}^{2p'}\Phi)dxdt+\int_{0}^{M}\int_{\R^{n}}\nabla u\nabla(\eta_{M}^{2p^{\prime}}\Phi)dxdt\\
&-\int_{0}^{M}\int_{\R^{n}}D_{1}(x)u\pa_{t}(\eta_{M}^{2p^{\prime}}\Phi)dxdt+\int^{M}_{0}\int_{\R^{n}}uV_{1}(x)\eta^{2p'}_{M}\Phi dxdt\\
=&\int_{0}^{M}\int_{\R^{n}}u\Big((\pa^{2}_{t}\eta^{2p'}_{M})\Phi+2(\pa_{t}\eta^{2p'}_{M})(\pa_{t}\Phi)-D_{1}(x)(\pa_{t}\eta^{2p'}_{M})\Phi\Big)dxdt\\
&+\int_{0}^{M}\int_{\R^{n}}u\eta^{2p'}_{M}\Big(\pa^{2}_{t}\Phi-\Delta \Phi-D_{1}(x)\pa_{t}\Phi+V_{1}(x)\Phi\big)\Big)dxdt\\
=&\int_{0}^{M}\int_{\R^{n}}u\Big((\pa^{2}_{t}\eta^{2p'}_{M})\Phi+2(\pa_{t}\eta^{2p'}_{M})(\pa_{t}\Phi)-D_{1}(x)(\pa_{t}\eta^{2p'}_{M})\Phi\Big)dxdt\\
=&I_{3}+I_{4}+I_{5}
\end{align*}
We apply H\"older inequality to the following estimates
\begin{align} \nonumber
|I_{3}|\lesssim&M^{-2}\int_{0}^{M}\int_{\R^{n}}|\theta^{2p'-2}_{M}\Phi|dxdt\\ \nonumber
\lesssim&M^{-2}\Big(\int_{0}^{M}\int_{\R^{n}}|u|^{q}\theta^{(2p'-2)q}_{M}dxdt\Big)^{\frac{1}{q}}\Big(\int_{\frac{M}{2}}^{M}\int_{\R^{n}}\Phi^{q'}dxdt\Big)^{\frac{1}{q'}}\\ \nonumber
\lesssim&M^{-2}\Big(\int_{0}^{M}\int_{\R^{n}}|u|^{q}\theta^{2p'}_{M}dxdt\Big)^{\frac{1}{q}}\Big(\int_{\frac{M}{2}}^{M}\int_{\R^{n}}\Phi^{q'}dxdt\Big)^{\frac{1}{q'}} \   (p<q)\\ \nonumber
\lesssim&M^{-2}\Big(\int_{0}^{M}\int_{\R^{n}}|u|^{q}\theta^{2p'}_{M}dxdt\Big)^{\frac{1}{q}}\Big(\int_{\frac{M}{2}}^{M}\int_{0}^{t+R}(r+1)^{n-1-\frac{n-1}{2}q'}e^{-(t-r)q'}drdt\Big)^{\frac{1}{q'}}\    (Lemma\ \ref{lp5})\\ \nonumber
\lesssim&M^{-2}\Big(\int_{0}^{M}\int_{\R^{n}}|u|^{q}\theta^{2p'}_{M}dxdt\Big)^{\frac{1}{q}}\Big(\int_{\frac{M}{2}}^{M}(2+t)^{n-1-\frac{n-1}{2}q'}dt\Big)^{\frac{1}{q'}}\  (Lemma\ \ref{le1})\\
\lesssim&M^{-2+(n-\frac{(n-1)q'}{2})\frac{1}{q'}}\Big(\int_{0}^{M}\int_{\R^{n}}|u|^{q}\theta^{2p'}_{M}dxdt\Big)^{\frac{1}{q}}, \nonumber
\end{align}
\begin{align}
|I_{4}|\lesssim M^{-1+(n-\frac{(n-1)q'}{2})\frac{1}{q'}}\Big(\int_{0}^{M}\int_{\R^{n}}|u|^{q}\theta^{2p'}_{M}dxdt\Big)^{\frac{1}{q}},\nonumber
\end{align}
\begin{align} \nonumber
|I_{5}|\lesssim&M^{-1}\int_{0}^{M}\int_{\R^{n}}|u\theta^{2p'-1}_{M}\Phi D_{1}(x)|dxdt\\ \nonumber
\lesssim&M^{-1}\Big(\int_{0}^{M}\int_{\R^{n}}|u|^{q}\theta^{2p'}_{M}dxdt\Big)^{\frac{1}{q}}\Big(\int_{\frac{M}{2}}^{M}\int_{\R^{n}}D^{q'}_{1}(x)\Phi^{q'}dxdt\Big)^{\frac{1}{q'}}\\ \nonumber
\lesssim&M^{-1}\Big(\int_{0}^{M}\int_{\R^{n}}|u|^{q}\theta^{2p'}_{M}dxdt\Big)^{\frac{1}{q}}\Big(\int_{\frac{M}{2}}^{M}(1+t)^{n-1-(\frac{n-1}{2}+\be)q'}dt\Big)^{\frac{1}{q'}} \ (Lemma\  \ref{le1})\\ 
\lesssim&M^{-2+(n-\frac{(n-1)q'}{2})\frac{1}{q'}}\Big(\int_{0}^{M}\int_{\R^{n}}|u|^{q}\theta^{2p'}_{M}dxdt\Big)^{\frac{1}{q}}. \nonumber
\end{align} 
This in turn implies that 
\begin{align}
\label{l8}
[C(u_{0},u_{1})\ep]^{q}\lesssim M^{\frac{(n-1)q}{2}-n}\Big(\int_{0}^{M}\int_{\R^{n}}|u|^{q}\theta^{2p'}_{M}dxdt\Big).
\end{align}
In an analogy way, we choose test function $\Psi=\eta^{2p'}_{M}\Phi(t,x)$ to replace $\Psi(t,x)$ in \eqref{l2}, we have that
\begin{align}
\label{l9}
[C(v_{0},v_{1})\ep]^{p}\lesssim M^{\frac{(n-1)p}{2}-n}\Big(\int_{0}^{M}\int_{\R^{n}}|v|^{p}\theta^{2p'}_{M}dxdt\Big) .
\end{align}
\eqref{l8} and \eqref{l7} yield 
\begin{align*}
M\lesssim\ep^{-\frac{1}{(p+q^{-1}+2)(qp-1)^{-1}-\frac{n-1}{2}}}.
\end{align*}
\eqref{l6} and \eqref{l9} yield 
\begin{equation*}
M\lesssim\ep^{-\frac{1}{(q+p^{-1}+2)(qp-1)^{-1}-\frac{n-1}{2}}}.
\end{equation*}
Note that, by replacing $p$ with $q$ in the test functions, we could have the corresponding estimates under the assumption $p\geq q$.  
Hence for any $M\in (1,T_{\ep})$, we get that
\begin{equation*}
M\lesssim \min\big\{\ep^{-\frac{1}{(p+q^{-1}+2)(qp-1)^{-1}-\frac{n-1}{2}}}, \ep^{-\frac{1}{(q+p^{-1}+2)(qp-1)^{-1}-\frac{n-1}{2}}}\big\}.
\end{equation*}  
That is to say
\begin{equation*}
T_{\ep}\lesssim \min\big\{\ep^{-\frac{1}{(p+q^{-1}+2)(qp-1)^{-1}-\frac{n-1}{2}}}, \ep^{-\frac{1}{(q+p^{-1}+2)(qp-1)^{-1}-\frac{n-1}{2}}}\big\}.
\end{equation*} 
So we come to the estimate
\begin{align*}
T_{\ep}\lesssim\ep^{-\Ga^{-1}_{SS}(n,p,q)}, \ \ \ \Ga_{SS}(n,p,q)>0.
\end{align*}
We obtain the first lifespan estimate in Theorem \ref{liu1}.
\subsubsection{{\bf Blow up region $\Ga_{SS}(n,p,q)=0$}}
Let
\begin{align}
\ \theta(t)=
\begin{cases}
\ 0,        \     &0\leq t<\frac{1}{2},\\
\eta(t),    &t\geq\frac{1}{2},
\end{cases}
\end{align}
where $\theta_{M}(t)=\theta(\frac{t}{M})$ and $M\in(1,T_{\ep})$. For nonnegative function $\omega \in L^{1}_{loc}([0,T); L^{1}(\R^{n}))$, we set
\begin{align*}
Y[\omega(t,x)](M)=\int^{M}_{1}\Big(\int^{T}_{0}\int_{\R^{n}}\omega(t,x)\theta^{2p'}_{\sigma}(t)dxdt\Big)\sigma^{-1}d\sigma.
\end{align*}
Then we can get that
\begin{equation}
\label{lp1}
Y'[\omega(t,x)](M)=\frac{d}{dM}Y[\omega(t,x)](M)=M^{-1}\Big(\int^{T}_{0}\int_{\R^{n}}\omega(t,x)\theta^{2p'}_{M}(t)dxdt\Big).
\end{equation}
By direct computation, we get that
\begin{align*}
Y[\omega(x)]&=\int_{1}^{M}\Big(\int_{0}^{T}\int_{\R^{n}}\omega(x)\theta_{\sigma}^{2p'}(t)dxdt\Big)\sigma^{-1}d\sigma=\int_{0}^{T}\int_{\R^{n}}\omega(x)\int_{\frac{t}{M}}^{t}\theta^{2p'}(s)s^{-1}dsdxdt\\
\lesssim&\int_{0}^{\frac{M}{2}}\int_{\R^{n}}\omega(x)\int_{\frac{t}{M}}^{t}\theta^{2p'}(s)s^{-1}dsdxdt+\int_{\frac{M}{2}}^{T}\int_{\R^{n}}\omega(x)\int_{\frac{t}{M}}^{t}\theta^{2p'}(s)s^{-1}dsdxdt\\
\lesssim&\int_{0}^{\frac{M}{2}}\int_{\R^{n}}\omega(x)\eta^{2p'}(\frac{t}{M})\int_{\frac{1}{2}}^{1}s^{-1}dsdxdt+\int_{\frac{M}{2}}^{T}\int_{\R^{n}}\omega(x)\theta^{2p'}(\frac{t}{M})
\int_{\frac{1}{2}}^{1}s^{-1}dsdxdt\\
\lesssim&\ln2\int^{T}_{0}\int_{\R^{n}}\omega(t,x)\eta^{2p'}_{M}(t)dxdt.
\end{align*}
Then we can obtain that
\begin{equation}
\label{lp2}
Y[\omega(t,x)](M)\lesssim\ln2\int^{T}_{0}\int_{\R^{n}}\omega(t,x)\eta^{2p'}_{M}(t)dxdt.
\end{equation}

{\bf {(\uppercase\expandafter{\romannumeral3}) Test function 3: $\Psi(t,x)=\eta^{2p'}_{M}(t)b_{a}(t,x)$}}
\\

Let $a=\frac{n-1}{2}-\frac{1}{q}>0$, by the lower bound of $b_{a}$ in Lemma \ref{le3}, we have that 
$$\int^{M}_0\int_{\R^{n}}|v|^{p}\theta^{2q'}_{M}b_{a}dxdt  \ \gt \ M^{-a} \int^{M}_0\int_{\R^{n}}|v|^{p}\theta^{2q'}_{M}dxdt\ .$$
By applying \eqref{ll} and \eqref{l8} (by replacing $p$ with $q$ in $\theta^{2p'}_M$), we know that
$$\int^{M}_0\int_{\R^{n}}|v|^{p}\theta^{2q'}_{M}dxdt \gt\ \ep^{pq}M^{n+p+1-\frac{n-1}{2}pq}\ .$$
Noting that when $p> q$, $\Ga_{SS}=(p+2+1/q)/(pq-1)$, then we have 
$$n+p+1-\frac{n-1}{2}pq-a=\Ga_{SS}(p, q, n)=0\ .$$
Hence we get that
\begin{equation}
\label{ll1}
\int^{M}_0\int_{\R^{n}}|v|^{p}\theta^{2q'}_{M}dxdt \gt\ \ep^{qp}\ .
\end{equation}
Next we choose test function $ \Psi=\eta^{2q'}_{M}b_{a}$ to replace $\Psi(t,x)$ in \eqref{l1} and applying \eqref{lll1}, we get that
\begin{align*} 
&\ep C(u_{0},u_{1})+\int_{0}^{M}\int_{\R^{n}}|v|^{p}\eta_{M}^{2q^{\prime}}(t)b_{a}dxdt\\
=&\int_{0}^{T}\int_{\R^{n}}u\Big((\pa^{2}_{t}\eta^{2q'}_{M})b_{a}+2(\pa_{t}\eta^{2q'}_{M})(\pa_{t}b_{a})-D_{1}(x)(\pa_{t}\eta^{2q'}_{M})b_{a}\Big)dxdt\\ 
=&I_{6}+I_{7}+I_{8}
\end{align*} 
We apply H\"older inequality to the following estimates
\begin{align} \nonumber
\ \ \ \ |I_{6}|\lesssim&M^{-2}\int_{0}^{M}\int_{\R^{n}}|u(t,x)\theta^{2q'-2}_{M}(t)b_{a}| dxdt\\ \nonumber
\lesssim&M^{-2}\Big(\int_{0}^{M}\int_{\R^{n}}|u|^{q}\theta^{2q'}_{M}dxdt\Big)^{\frac{1}{q}}\Big(\int_{\frac{M}{2}}^{M}\int_{\R^{n}}b_{a}^{q'}dxdt\Big)^{\frac{1}{q'}} \ \ \\ \nonumber
\lesssim&M^{-2}\Big(\int_{0}^{M}\int_{\R^{n}}|u|^{q}\theta^{2q'}_{M}dxdt\Big)^{\frac{1}{q}}\Big(\int_{\frac{M}{2}}^{M}\int_{0}^{t+R}(t+R)^{(\frac{1}{q}-\frac{n-1}{2})q'}r^{n-1}drdt\Big)^{\frac{1}{q'}} \ \ \ \ \ (Lemma \ \ref{le3})\\ \nonumber
\lesssim&M^{\frac{nq-2n-q}{2q}}\Big(\int_{0}^{M}\int_{\R^{n}}|u|^{q}\theta^{2q'}_{M}dxdt\Big)^{\frac{1}{q}}, 
\end{align}

\begin{align}\nonumber
|I_{7}|\lesssim&M^{-1}\Big(\int_{0}^{M}\int_{\R^{n}}|u|^{q}\theta^{2q'}_{M}dxdt\Big)^{\frac{1}{q}}\Big(\int_{\frac{M}{2}}^{M}\int_{\R^{n}}b_{a+1}^{q'}dxdt\Big)^{\frac{1}{q'}} \\ \nonumber
\lesssim&M^{\frac{nq-2n-q}{2q}}(\ln M)^{\frac{q-1}{q}}\Big(\int_{0}^{M}\int_{\R^{n}}|u|^{q}\theta^{2q'}_{M}dxdt\Big)^{\frac{1}{q}}, \nonumber
\end{align} 
\begin{align} \nonumber
\ \ \ \ |I_{8}|\lesssim&M^{-1}\int_{0}^{M}\int_{\R^{n}}|u(t,x)\theta^{2q'-1}_{M}(t)b_{a}D_{1}(x)|dxdt\\ \nonumber 
\lesssim&M^{-1}\Big(\int_{0}^{M}\int_{\R^{n}}|u|^{q}\theta^{2q'}_{M}(t)dxdt\Big)^{\frac{1}{q}}\Big(\int_{\frac{M}{2}}^{M}\int_{\R^{n}}D^{q'}_{1}(x)b_{a}^{q'}dxdt\Big)^{\frac{1}{q'}}\\ \nonumber 
\lesssim&M^{-1}\Big(\int_{0}^{M}\int_{\R^{n}}|u|^{q}\theta^{2q'}_{M}(t)dxdt\Big)^{\frac{1}{q}}\Big(\int_{\frac{M}{2}}^{M}\int^{t+R}_{0}(R+t)^{{(\frac{1}{q}-\frac{n-1}{2})q'}}r^{n-1}(1+r)^{-\be q'}drdt\Big)^{\frac{1}{q'}}\\  \nonumber
\lesssim&M^{\frac{nq+q-2n-2\be q}{2q}}\Big(\int_{0}^{M}\int_{\R^{n}}|u|^{q}\theta^{2q'}_{M}(t)dxdt\Big)^{\frac{1}{q}} \\ \nonumber
\lesssim&M^{\frac{nq-2n-q}{2q}}\Big(\int_{0}^{M}\int_{\R^{n}}|u|^{q}\theta^{2q'}_{M}(t)dxdt\Big)^{\frac{1}{q}}. \nonumber
\end{align}
Hence, we get that
\begin{align}
\label{l10}
\int_{0}^{M}\int_{\R^{n}}|v|^{p}\eta_{M}^{2q^{\prime}}(t)b_{a}dxdt\lesssim M^{\frac{nq-2n-q}{2q}}(\ln M)^{\frac{q-1}{q}}\Big(\int_{0}^{M}\int_{\R^{n}}|u|^{q}\theta^{2q'}_{M}dxdt\Big)^{\frac{1}{q}},
\end{align}

\begin{equation}
\label{l11}
\int_{0}^{M}\int_{\R^{n}}|u|^{q}\eta_{M}^{2q^{\prime}}(t)b_{a}dxdt\lesssim M^{\frac{np-2n-p}{2p}}(\ln M)^{\frac{p-1}{p}}\Big(\int_{0}^{M}\int_{\R^{n}}|v|^{p}\theta^{2q'}_{M}dxdt\Big)^{\frac{1}{p}}.
\end{equation}
Plugging \eqref{ll} into \eqref{l10} and by the lower bound of $b_{a}$ in Lemma \ref{le3}, we get that
\begin{equation}
\label{ll2}
(\int^{M}_{0}\int_{\R^{n}}|v|^{p}\eta^{2q'}_{M}b_{a}dxdt)^{pq}\lesssim (\ln M)^{p(q-1)}\int^{M}_{0}\int_{\R^{n}}|v|^{p}\theta^{2q'}_{M}b_{a}dxdt .
\end{equation}
Let 
$$Y[b_{a}|v|^{p}](M)=\int^{M}_{1}(\int^{T}_{0}\int_{\R^{n}}|v|^{p}b_{a}\theta^{2q'}_{\sigma}(t)dxdt)\sigma^{-1}d\sigma$$
By combining \eqref{lp1}, \eqref{lp2}, \eqref{ll1} and \eqref{ll2}, we have that
\begin{align*}
\ep^{qp}\lesssim& MY'[b_{a}|v|^{p}](M)\\
[Y[b_{a}|v|^{p}](M)]^{pq}\lesssim& (\ln M)^{p(q-1)}MY'[b_{a}|v|^{p}](M)
\end{align*}
Exploiting Lemma \ref{le4} with $ p_{2}=p(q-1)+1,$ $p_{1}=pq,$ $ \delta=\ep^{pq},$ we have that
\begin{align*}
M\lesssim\exp(\ep^{-q(pq-1)}),\ \ p>q.
\end{align*}
Hence for any $M \in (1,T_{\ep}),$ that is to say
\begin{align*}
T_{\ep}\lesssim\exp(\ep^{-q(pq-1)}),\ \ p>q.
\end{align*}
By the symmetrical characteristic of $p$ and $q$, we could get that 
 \begin{align*}
T_{\ep}\lesssim\exp(\ep^{-p(pq-1)}), \ \ p<q.
\end{align*}
Then we get that
\begin{equation*}
T_{\ep}\les\exp\big(\ep^{-\min\{p(pq-1),q(pq-1)\} }\big),\  \ \Ga_{SS}(n,p,q)=0,\ p\neq q.
\end{equation*}
Next, we give the lifespan estimate when $\Ga_{SS}(n,p,q)=0$, $p=q=p_{S}(n)$. In the following, we take $D(x)=D_1(x)+D_2(x)$. By combining \eqref{l8}, \eqref{l9}, we have
\begin{align*}
\ep &C(u_{0},u_{1})+\ep C(v_{0},v_{1})+\int^{M}_{0}\int_{\R^{n}}|u+v|^{p}\eta_{M}^{2p'}\Phi dxdt\\ 
\lesssim\ &\ep C(u_{0},u_{1})+ \ep C(v_{0},v_{1})+\int^{M}_{0}\int_{\R^{n}}(|u|^{p}+|v|^{p})\eta_{M}^{2p'}\Phi dxdt\\ 
\lesssim\ &\int_{0}^{M}\int_{\R^{n}}(u+v)\Big((\pa^{2}_{t}\eta^{2p'}_{M})\Phi+2(\pa_{t}\eta^{2p'}_{M})(\pa_{t}\Phi)-D(x)(\pa_{t}\eta^{2p'}_{M})\Phi\Big)dxdt\ \\
\lesssim\ &M^{\frac{(n-1)}{2}-\frac{n}{p}}\Big(\int^{M}_{0}\int_{\R^{n}}|u+v|^{p}\theta^{2p'}_{M}dxdt\Big)^{\frac{1}{p}}\\
\lesssim \ &M^{\frac{(n-1)}{2}-\frac{n}{p}}M^{\frac{n-1}{2p}-\frac{1}{p^{2}}}\big(\int^{M}_{0}\int_{\R^{n}}|u+v|^{p}\theta^{2p'}_{M}b_{a}dxdt\big)^{\frac{1}{p}}\\
\end{align*}
Noting that when $\Ga_{SS}(n,p,q)=0$, $p=q=p_{S}(n)$, we have that
\begin{equation*}
n-\frac{(n-1)p}{2}=\frac{n-1}{2}-\frac{1}{p}.
\end{equation*}
Which leads to
\begin{align}
\label{f1}
\ep^{p}\lesssim\big(\int^{M}_{0}\int_{\R^{n}}|u+v|^{p}\theta^{2p'}_{M}b_{a}dxdt\big).
\end{align}
Similarly, we can get that
\begin{align*}
&\ep +\int_{0}^{M}\int_{\R^{n}}|v+u|^{p}\eta_{M}^{2p^{\prime}}(t)b_{a}dxdt\\
\lesssim\ &\int_{0}^{M}\int_{\R^{n}}(u+v)\Big((\pa^{2}_{t}\eta^{2p'}_{M})b_{a}+2(\pa_{t}\eta^{2p'}_{M})(\pa_{t}b_{a})-D(x)(\pa_{t}\eta^{2p'}_{M})b_{a}\Big)dxdt\ \\ 
\lesssim\ &I_{9}+I_{10}+I_{11}
\end{align*}
We apply H\"older inequality to the following estimates
\begin{align*}
|I_{9}|\lesssim&M^{-2}\Big(\int^{M}_{0}\int_{\R^{n}}|u+v|^{p}\theta^{2p'}_{M}b_{a}dxdt\Big)^{\frac{1}{p}}\Big(\int^{M}_{\frac{M}{2}}\int_{\R^{n}}b_{a}dxdt\Big)^{\frac{p-1}{p}}\\
\lesssim&M^{-2}\Big(\int^{M}_{\frac{M}{2}}\int^{t+R}_{0}(t+R)^{-(\frac{n-1}{2}-\frac{1}{p})}r^{n-1}drdt\Big)^{\frac{p-1}{p}}\Big(\int^{M}_{0}\int_{\R^{n}}|u+v|^{p}\theta^{2p'}_{M}b_{a}dxdt\Big)^{\frac{1}{p}}\\
\lesssim&\Big(\int^{M}_{0}\int_{\R^{n}}|u+v|^{p}\theta^{2p'}_{M}b_{a}dxdt\Big)^{\frac{1}{p}},
\end{align*}
\begin{align*}
|I_{10}|\lesssim&M^{-1}\Big(\int^{M}_{0}\int_{\R^{n}}|u+v|^{p}\theta^{2p'}_{M}b_{a}dxdt\Big)^{\frac{1}{p}}\Big(\int^{M}_{\frac{M}{2}}\int_{\R^{n}}b^{-(\frac{1}{p-1})}_{a}b^{\frac{p}{p-1}}_{a}dxdt\Big)^{\frac{p-1}{p}}\\
\lesssim&M^{-1}\Big(\int^{M}_{\frac{M}{2}}\int^{t+R}_{0}(t+R)^{(\frac{n-1}{2}-\frac{1}{p})\frac{1}{p-1}}(t+R)^{-(\frac{n-1}{2})\frac{p}{p-1}}(t+R+1-r)^{-1}r^{n-1}drdt\Big)^{\frac{p-1}{p}}\\
&\times\Big(\int^{M}_{0}\int_{\R^{n}}|u+v|^{p}\theta^{2p'}_{M}b_{a}dxdt\Big)^{\frac{1}{p}}\\
\lesssim&M^{-1}\Big(\int^{M}_{\frac{M}{2}}(t+R)^{\frac{n-1}{2}-\frac{1}{p(p-1)}}\ln (t+R)\Big)^{\frac{p-1}{p}}\Big(\int^{M}_{0}\int_{\R^{n}}|u+v|^{p}\theta^{2p'}_{M}b_{a}dxdt\Big)^{\frac{1}{p}}\\
\lesssim&(\ln M)^{\frac{p-1}{p}}\Big(\int^{M}_{0}\int_{\R^{n}}|u+v|^{p}\theta^{2p'}_{M}b_{a}dxdt\Big)^{\frac{1}{p}},
\end{align*}
\begin{align*}
|I_{11}|\lesssim&M^{-1}\Big(\int^{M}_{0}\int_{\R^{n}}|u+v|^{p}\theta^{2p'}_{M}b_{a}dxdt\Big)^{\frac{1}{p}}\Big(\int^{M}_{\frac{M}{2}}\int_{\R^{n}}D^{p'}(x)b_{a}dxdt\Big)^{\frac{p-1}{p}}\\
\lesssim&M^{-1}\Big(\int^{M}_{\frac{M}{2}}\int^{t+R}_{0}(t+R)^{-(\frac{n-1}{2}-\frac{1}{p})}(1+r)^{-\be p'}r^{n-1}drdt\Big)^{\frac{p-1}{p}}\Big(\int^{M}_{0}\int_{\R^{n}}|u+v|^{p}\theta^{2p'}_{M}b_{a}dxdt\Big)^{\frac{1}{p}}\\
\lesssim&\Big(\int^{M}_{0}\int_{\R^{n}}|u+v|^{p}\theta^{2p'}_{M}b_{a}dxdt\Big)^{\frac{1}{p}}.
\end{align*}
Therefor we conclude from the estimates $I_{9}$-$I_{11}$ that
\begin{equation}
\label{f2}
\Big(\int_{0}^{M}\int_{\R^{n}}|v+u|^{p}\eta_{M}^{2p^{\prime}}(t)b_{a}dxdt\Big)^{p}\lesssim
(\ln M)^{p-1}\Big(\int^{M}_{0}\int_{\R^{n}}|u+v|^{p}\theta^{2p'}_{M}b_{a}dxdt\Big).
\end{equation}
Let 
$$Y[b_{a}|u+v|^{p}](M)=\int^{M}_{1}\Big(\int^{M}_{0}\int_{\R^{n}}|u+v|^{p}b_{a}\theta^{2q'}_{\sigma}(t)dxdt\Big)\sigma^{-1}d\sigma.$$
By combining \eqref{lp1}, \eqref{lp2}, \eqref{f1} and \eqref{f2}, we get that
\begin{align*}
\ep^{p}\lesssim& MY'[b_{a}|u+v|^{p}](M)\\
[Y[b_{a}|u+v|^{p}](M)]^{p}\lesssim& (\ln M)^{p-1}MY'[b_{a}|u+v|^{p}](M)
\end{align*}
By exploiting Lemma \ref{le4} with $p_{2}=p_{1}=p,$ $\delta=\ep^{p},$ we have that
\begin{align*}
M\lesssim\exp(\ep^{-p(p-1)}), \ \ p=q=p_{S}(n).
\end{align*}
Hence for any $M \in (1,T_{\ep}),$ that is to say
\begin{align*}
T_{\ep}\lesssim\exp(\ep^{-p(p-1)}), \ \ p=q=p_{S}(n).
\end{align*}
We finish the proof of Theorem \ref{liu1}.

\subsection{Proof of Theorem \ref{liu2}}
\subsubsection{{\bf Blow up region\ $\Ga_{GG}(n,p,q)>0$}}\ \\

{\bf {(\uppercase\expandafter{\romannumeral1}) Test function 1: $\Psi(t,x)=\eta^{2p'}_{M}(t)\phi_0(x)$}}
\\

Firstly, we choose test function $\Psi=\eta^{2p'}_{M}(t)\phi_0$ to replace $\Psi(t,x)$ in \eqref{l1} we get
\begin{align}\nonumber
&\ep\int_{\R^{n}}u_{1}(x)\phi_0 dx+\ep\int_{\R^{n}}u_{0}(x)D_{1}(x)\phi_0 dx+\int_{0}^{M}\int_{\R^{n}}|v_{t}|^{p}\eta_{M}^{2p^{\prime}}(t)\phi_0 dxdt\\  \nonumber
=&-\int_{0}^{M}\int_{\R^{n}}u_{t}\phi_0(\pa_{t} \eta_{M}^{2p'}(t))dxdt+\int_{0}^{M}\int_{\R^{n}}u_{t}(D_{1}(x)) \eta_{M}^{2p'}(t)\phi_0dxdt\\
&-\int_{0}^{M}\int_{\R^{n}}u\eta_{M}^{2p'}(\Delta \phi_0-V_{1}(x)\phi_0)dxdt \\ \nonumber
=&I_{12}+I_{13}
\end{align}
We apply H\"older inequality to the following estimates
\begin{align}\nonumber
|I_{12}|\lesssim&M^{-1}\int^{M}_{0}\int_{\R^{n}}|u_{t}\eta^{2p'-1}_{M}|dxdt\\
\lesssim&M^{\frac{nq-n-1}{q}}\Big(\int^{M}_{0}\int_{\R^{n}}|u_{t}|^{q}\eta^{2p'}_{M}dxdt\Big)^{\frac{1}{q}},\nonumber
\end{align}
\begin{align} \nonumber
|I_{13}|\lesssim&\int^{M}_{0}\int_{\R^{n}}|u_{t}\eta^{2p'}_{M}D_{1}(x)|dxdt\\ \nonumber
\lesssim&\Big(\int^{M}_{0}\int_{\R^{n}}|u_{t}|^{q}\eta^{2p'}_{M}dxdt\Big)^{\frac{1}{q}}\Big(\int^{M}_{\frac{M}{2}}\int_{\R^{n}}(D_{1}(x))^{q'}dxdt\Big)^{\frac{1}{q'}}\\ \nonumber
\lesssim&\Big(\int^{M}_{0}\int_{\R^{n}}|u_{t}|^{q}\eta^{2p'}_{M}dxdt\Big)^{\frac{1}{q}}
\begin{cases}
\ (M\ln(M))^{\frac{1}{q'}}, \ \ \  &n-\be q'=0,\\
\ M^{(n+1-\be q')\frac{1}{q'}},\ \ \  &n-\be q'>0,\\
\ M^{\frac{1}{q'}},\ \ \  &n-\be q'<0.
\end{cases}\\
\lesssim&M^{\frac{nq-n-1}{q}}\Big(\int^{M}_{0}\int_{\R^{n}}|u_{t}|^{q}\eta^{2p'}_{M}dxdt\Big)^{\frac{1}{q}} .\nonumber
\end{align}
So we get that
\begin{equation}
\label{ll12}
\int_{0}^{M}\int_{\R^{n}}|v_{t}|^{p}\eta_{M}^{2p^{\prime}}dxdt\lesssim M^{\frac{nq-n-1}{q}}\Big(\int^{M}_{0}\int_{\R^{n}}|u_{t}|^{q}\eta^{2p'}_{M}dxdt\Big)^{\frac{1}{q}} ,
\end{equation}
Similarly, we have that
\begin{equation}
\label{l13}
\int_{0}^{M}\int_{\R^{n}}|u_{t}|^{q}\eta_{M}^{2p^{\prime}}dxdt\lesssim M^{\frac{np-n-1}{p}}\Big(\int^{M}_{0}\int_{\R^{n}}|v_{t}|^{p}\eta^{2p'}_{M}dxdt\Big)^{\frac{1}{p}} .
\end{equation}
Applying \eqref{ll12} and \eqref{l13} yields
\begin{align}
\int_{0}^{M}\int_{\R^{n}}|v_{t}|^{p}\eta_{M}^{2p^{\prime}}dxdt\lesssim M^{\frac{nqp-p-n-1}{qp-1}},\\
\int_{0}^{M}\int_{\R^{n}}|u_{t}|^{q}\eta_{M}^{2p^{\prime}}dxdt\lesssim M^{\frac{nqp-q-n-1}{qp-1}}.
\end{align}

{\bf {(\uppercase\expandafter{\romannumeral2}) Test function 2: $\Psi=-\pa_{t}(\eta^{2p'}_{M}(t)\Phi(t,x))$}}

Secondly, we choose another test function $ \Psi=-\pa_{t}(\eta^{2p'}_{M}\Phi)$ to replace $\Psi$ in \eqref{l1}.
It is worth observing that
\begin{align*}
\Psi(t,x)=-\pa_{t}(\eta^{2p'}_{M}\Phi)=\eta^{2p'}_{M}\Phi-2p'\eta^{2p'-1}_{M}(\pa_{t}\eta_{M})\Phi\geq \eta^{2p'}_{M}\Phi>0\ .
\end{align*}
Direct calculation shows  
\begin{align} \nonumber
\ep C(u_{0}&,u_{1})+\int_{0}^{M}\int_{\R^{n}}|v_{t}|^{p}(\eta_{M}^{2p^{\prime}}\Phi-2p'\eta^{2p'-1}_{M}(\pa_{t}\eta_{M})\Phi)dxdt\\  \nonumber
&=\int^{M}_{0}\int_{\R^{n}}u_{t}\Big((\pa^{2}_{t}\eta^{2p'}_{M})\Phi+2(\pa_{t}\eta^{2p'}_{M})(\pa_{t}\Phi)-D_{1}(x)(\pa_{t}\eta^{2p'}_{M})\Phi\Big)dxdt\\
&+\int_{0}^{M}\int_{\R^{n}}u_t\eta^{2p'}_{M}\Big(\pa^{2}_{t}\Phi-\Delta \Phi-D_{1}(x)\pa_{t}\Phi+V_{1}(x)\Phi\big)\Big)dxdt\\
&= I_{14}+I_{15}+I_{16}\nonumber
\end{align}
By the same procedure in \eqref{l3}-\eqref{l4}, we get that
\begin{align}
\label{l14}
\ep C(u_{0},u_{1})\lesssim M^{-1+(n-(\frac{n-1}{2})q')\frac{1}{q'}}\Big(\int^{M}_{0}\int_{\R^{n}}|u_{t}|^{q}\eta^{2p'}_{M}dxdt\Big)^{\frac{1}{q}},
\end{align}
\begin{align}
\label{l15}
\ep C(v_{0},v_{1})\lesssim M^{-1+(n-(\frac{n-1}{2})p')\frac{1}{p'}}\Big(\int^{M}_{0}\int_{\R^{n}}|v_{t}|^{q}\eta^{2p'}_{M}dxdt\Big)^{\frac{1}{p}}.
\end{align}
On the one hand, by \eqref{ll12}, \eqref{l13}, \eqref{l14} and \eqref{l15}, we can obtain the  first lifespan estimation  in Theorem \ref{liu2}. 

On the other hand, by exploiting the estimates of $\Phi$, we have that
\begin{align*}
|I_{14}|\lesssim&\int_{0}^{M}\int_{\R^{n}}|u_{t}(\pa^{2}_{t}\eta^{2p'}_{M})\Phi|dxdt\\
\lesssim&M^{-2}\Big(\int_{0}^{M}\int_{\R^{n}}\theta^{2p'}_{M}|u_{t}|^{q}\Phi dxdt\Big)^{\frac{1}{q}}\Big(\int^{M}_{\frac{M}{2}}\int_{\R^{n}}\Phi dxdt\Big)^{\frac{1}{q'}}\\
\lesssim&M^{-2}\Big(\int_{0}^{M}\int_{\R^{n}}\theta^{2p'}_{M}|u_{t}|^{q}\Phi dxdt\Big)^{\frac{1}{q}}\Big(\int_{\frac{M}{2}}^{M}(t+1)^{\frac{n-1}{2}}dxdt\Big)^{\frac{1}{q'}}\\
\lesssim&M^{\frac{nq-n-3q-1}{2q}}\Big(\int_{0}^{M}\int_{\R^{n}}\theta^{2p'}_{M}|u_{t}|^{q}\Phi dxdt\Big)^{\frac{1}{q}},
\end{align*}
\begin{align*}
|I_{15}|\lesssim&\int_{0}^{M}\int_{\R^{n}}|u_{t}(\pa_{t}\eta^{2p'}_{M})\Phi|dxdt\\
\lesssim&M^{-1}\Big(\int_{0}^{M}\int_{\R^{n}}\theta^{2p'}_{M}|u_{t}|^{q}\Phi dxdt\Big)^{\frac{1}{q}}\Big(\int^{M}_{\frac{M}{2}}\int_{\R^{n}}\Phi dxdt\Big)^{\frac{1}{q'}}\\
\lesssim&M^{\frac{nq-n-q-1}{2q}}\Big(\int_{0}^{M}\int_{\R^{n}}\theta^{2p'}_{M}|u_{t}|^{q}\Phi dxdt\Big)^{\frac{1}{q}},
\end{align*}
\begin{align*}
|I_{16}|\lesssim|I_{15}|\lesssim M^{\frac{nq-n-q-1}{2q}}\Big(\int_{0}^{M}\int_{\R^{n}}\theta^{2p'}_{M}|u_{t}|^{q}\Phi dxdt\Big)^{\frac{1}{q}}.
\end{align*}
Then we get that
\begin{align}
\label{lml1}
\ep C(u_{0}&,u_{1})+\int_{0}^{M}\int_{\R^{n}}|v_{t}|^{p}\eta_{M}^{2p^{\prime}}\Phi dxdt\lesssim M^{\frac{nq-n-1-q}{2q}}\Big(\int_{0}^{M}\int_{\R^{n}}|u_{t}|^{q}\theta_{M}^{2p^{\prime}}\Phi dxdt\Big)^{\frac{1}{q}},
\end{align}
\begin{align}
\label{lml2}
\ep C(v_{0}&,v_{1})+\int_{0}^{M}\int_{\R^{n}}|u_{t}|^{q}\eta_{M}^{2p^{\prime}}\Phi dxdt\lesssim M^{\frac{np-n-1-p}{2p}}\Big(\int_{0}^{M}\int_{\R^{n}}|v_{t}|^{p}\theta_{M}^{2p^{\prime}}\Phi dxdt\Big)^{\frac{1}{p}}.
\end{align}
Combining above two inequalities, we have that
\begin{align}
\label{ll4}
\ep C(u_{0}&,u_{1})+\int_{0}^{M}\int_{\R^{n}}|v_{t}|^{p}\eta_{M}^{2p^{\prime}}\Phi dxdt\lesssim M^{\frac{nqp-2p-n-qp-1}{2qp}}\Big(\int_{0}^{M}\int_{\R^{n}}|v_{t}|^{p}\theta_{M}^{2p^{\prime}}\Phi dxdt\Big)^{\frac{1}{qp}},\\
\ep C(v_{0}&,v_{1})+\int_{0}^{M}\int_{\R^{n}}|u_{t}|^{q}\eta_{M}^{2p^{\prime}}\Phi dxdt\lesssim M^{\frac{np-2q-n-qp-1}{2qp}}\Big(\int_{0}^{M}\int_{\R^{n}}|u_{t}|^{q}\theta_{M}^{2p^{\prime}}\Phi dxdt\Big)^{\frac{1}{qp}}.
\end{align}
We set $ A=\int_{0}^{M}\int_{\R^{n}}|v_{t}|^{p}\eta_{M}^{2p^{\prime}}\Phi dxdt$,  $B=\int_{0}^{M}\int_{\R^{n}}|u_{t}|^{q}\eta_{M}^{2p^{\prime}}\Phi dxdt$, then by above two inequalities, we apply the Young inequality
\begin{align*}
\ep+A\lesssim M^{\frac{nqp-2p-n-qp-1}{2qp}}A^{\frac{1}{qp}}
\leq\frac{A}{pq}+CM^{\frac{nqp-2p-n-qp-1}{2qp}(pq)'}\frac{1}{(pq)'}.
\end{align*}
Then
\begin{align*}
\ep+(1-\frac{1}{pq})A\lesssim M^{\frac{nqp-2p-n-qp-1}{2qp}(pq)'}\frac{1}{(pq)'}\ ,
\end{align*}
which lead to
\begin{align*}
M\leq T_{\ep}\lesssim \ep^{-\frac{1}{\frac{p+1}{qp-1}-\frac{n-1}{2}}},\ \ p<q.
\end{align*}
By the symmetrical characteristic of $p$ and $q,$ we get that
\begin{align*}
T_{\ep}\lesssim
\begin{cases}
\ep^{-\frac{1}{(\frac{p+1}{qp-1}-\frac{n-1}{2})}}, \ \ \ p<q,\\
\ep^{-\frac{1}{(\frac{q+1}{qp-1}-\frac{n-1}{2})}}, \ \ \ q<p.
\end{cases}
\end{align*}
So we come to the estimate
\begin{align*}
T_{\ep}\lesssim\ep^{-\Ga^{-1}_{GG}(n,p,q)},\ \ \ \Ga_{GG}(n,p,q)>0.
\end{align*}
Which is first part of lifespan estimate in Theorem \ref{liu2}.
\subsubsection{{\bf Blow up region\ $\Ga_{GG}(n,p,q)=0$}}
For the case $\Ga_{GG}(n,p,q)=0$, $p\neq q$. By \eqref{ll4} and  \eqref{lp1}, \eqref{lp2}, we have
\begin{align}
[C\ep+Y[\Phi|v_{t}|^{p}](M)]^{qp}\lesssim MY'[\Phi|v_{t}|^{p}](M)\ .
\end{align}
Let $Z(M)=C\ep+Y[\Phi|v_{t}|^{p}](M)$, then $Z(1)=C\ep$ and $Z'(M)=Y'(M)$, by above inequality, we have 
$$MZ'(M)\ \gt \ Z^{pq}\ , Z(1)=C\ep\ .$$
This is a typical ODE inequality which will blow up in finite time. By direct computation, we have that 
\begin{align*}
M\leq T_{\ep}\lesssim \exp(\ep^{-(qp-1)}).
\end{align*}
This is the second lifespan estimate in Theorem \ref{liu2}. For the case $\Ga_{GG}(n,p,q)=0$, $p=q=p_{G}(n)$, we set
\begin{align}
Y[\Phi|v_{t}+u_{t}|^{p}](M)=\int^{M}_{1}\Big(\int^{T}_{0}\int_{\R^{n}}|v_{t}+u_{t}|^{p}\Phi\theta^{2p'}_{\sigma}(t)dxdt\Big)\sigma^{-1}d\sigma.
\end{align}
From \eqref{lml1} and \eqref{lml2}, we get that
\begin{align}
[C\ep+Y[\Phi|v_{t}+u_{t}|^{p}]^{p}\lesssim MY'[\Phi|v_{t}+u_{t}|^{p}](M)\ .
\end{align}
Which lead to
\begin{align*}
M\leq T_{\ep}\lesssim \exp(\ep^{-(p-1)}),
\end{align*}
by similar procedure above. We have completed all the proofs of Theorem \ref{liu2}.

\subsection{Proof of Theorem \ref{liu3}}
\subsubsection{{\bf Blow up region\ $\Ga_{SG}(n,p,q)>0$}}\  
\\

{\bf {(\uppercase\expandafter{\romannumeral1})Test function 1: $\Psi(t,x)=\theta^{2p'}_{M}(t)\phi_0(x), \Psi(t, x)=\eta^{2p'}_{M}(t)\phi_0(x)$}}

Firstly, we choose test function $\Psi=\theta^{2p'}_{M}(t)\phi_0(x),$ $\eta^{2p'}_{M}(t)\phi_0(x)$. By substituting $ \Psi(t,x) $ in \eqref{l1} and \eqref{l2} respectively, and applying Lemma \ref{le5}, we get that
\begin{align*}
&\int_{0}^{M}\int_{\R^{n}}|v|^{q}\theta_{M}^{2p^{\prime}}\phi_0 dxdt\\  \nonumber
=&\int_{0}^{M}\int_{\R^{n}}u_{t}\phi_0\Big(-\pa_{t} \theta_{M}^{2p'}+D_{1}(x)\theta_{M}^{2p'}\Big)dxdt\\  \nonumber
=&I_{17}+I_{18}
\end{align*}
\begin{align*}
\ \ \ \ \ \ \ &\ep\int_{\R^{n}}v_{1}(x)\phi_0 dx+\ep\int_{\R^{n}}v_{0}(x)D_{2}(x)\phi_{0} dx+\int_{0}^{M}\int_{\R^{n}}|u_{t}|^{p}\eta_{M}^{2p^{\prime}}\phi_{0} dxdt\\  \nonumber
=&\int_{0}^{M}\int_{\R^{n}}v\phi_{0}\Big(\pa_{t}^{2}\eta_{M}^{2p'}-D_{2}(x)\pa_{t}\eta_{M}^{2p'}\Big)dxdt\\  \nonumber
=&I_{19}+I_{20}
\end{align*}
By H\"older's inequality, we have
\begin{align}
\label{l16}
\int_{0}^{M}\int_{\R^{n}}|v|^{q}\theta_{M}^{2p^{\prime}}dxdt\lesssim M^{\frac{np-n-1}{p}}\Big(\int_{0}^{M}\int_{\R^{n}}|u_{t}|^{p}\theta_{M}^{2p^{\prime}}dxdt\Big)^{\frac{1}{p}},
\end{align}
and
\begin{align}
\label{l17}
\int_{0}^{M}\int_{\R^{n}}|u_{t}|^{p}\eta_{M}^{2p^{\prime}}dxdt\lesssim M^{\frac{nq-n-q-1}{q}}\Big(\int_{0}^{M}\int_{\R^{n}}|v|^{q}\theta_{M}^{2p^{\prime}}dxdt\Big)^{\frac{1}{q}}.
\end{align}
Then by above two inequalities, we can get that
\begin{align}
\label{pp1}
\int_{0}^{M}\int_{\R^{n}}|v|^{q}\theta_{M}^{2p^{\prime}}dxdt\lesssim M^{\frac{npq-n-2q-1}{pq-1}},
\end{align}

\begin{align}
\label{pp2}
\int_{0}^{M}\int_{\R^{n}}|u_{t}|^{p}\eta_{M}^{2p^{\prime}}dxdt\lesssim M^{\frac{npq-n-pq-p-1}{pq-1}}.
\end{align}
{\bf {(\uppercase\expandafter{\romannumeral2})
Test function 2: $\Psi(t,x)=-\pa_{t}(\eta^{2p'}_{M}(x)\Phi(t,x)), \eta^{2p'}_{M}\Phi(t,x)$}}

Secondly, we choose another test function $\Psi=-\pa_{t}(\eta^{2p'}_{M}\Phi)$. By substituting $\Psi(t,x)$ in \eqref{l1} and $\Psi=\eta^{2p'}_{M}\Phi$ in \eqref{l2}, we get that
\begin{align*} \nonumber
\ep&C(u_{0},u_{1})+\int_{0}^{M}\int_{\R^{n}}|v|^{q}\pa_{t}(-\eta^{2p'}_{M}\Phi)dxdt\\  
&=\int^{M}_{0}\int_{\R^{n}}u_{t}\Big((\pa^{2}_{t}\eta^{2p'}_{M})\Phi+2(\pa_{t}\eta^{2p'}_{M})(\pa_{t}\Phi)-D_{1}(x)(\pa_{t}\eta^{2p'}_{M})\Phi\Big)dxdt\\ \nonumber
&= I_{21}+I_{22}+I_{23}
\end{align*}
\begin{align*} \nonumber
\ \ \ \ &\ep C(v_{0},v_{1})+\int_{0}^{M}\int_{\R^{n}}|u_{t}|^{p}\eta^{2p'}_{M}\Phi dxdt\\  
&=\int^{M}_{0}\int_{\R^{n}}v\Big((\pa^{2}_{t}\eta^{2p'}_{M})\Phi+2(\pa_{t}\eta^{2p'}_{M})(\pa_{t}\Phi)-D_{2}(x)(\pa_{t}\eta^{2p'}_{M})\Phi\Big)dxdt\\ \nonumber
&= I_{24}+I_{25}+I_{26}
\end{align*}
Then we apply H\"older inequality to the above terms, we can get that
\begin{align}
\label{l18}
[\ep C(u_{0},u_{1})]^{p}\lesssim M^{\frac{np-p-2n}{2}}\int_{0}^{M}\int_{\R^{n}}|u_{t}|^{p}\theta_{M}^{2p^{\prime}}dxdt,
\end{align}
\begin{align}
\label{l19}
[\ep C(v_{0},v_{1})]^{q}\lesssim M^{\frac{nq-q-2n}{2}}\int_{0}^{M}\int_{\R^{n}}|v|^{q}\theta_{M}^{2p^{\prime}}dxdt.
\end{align}
Hence we have by combining \eqref{pp2} and \eqref{l18} that
\begin{align*}
M\lesssim\ep^{-F^{-1}_{SG,1}(n,p,q)}.
\end{align*}
Utilizing \eqref{pp1} and \eqref{l19}, we obtain the estimate
\begin{align*}
M\lesssim\ep^{-F^{-1}_{SG,2}(n,p,q)}.
\end{align*}
Since for any $M \in (1,T_{\ep})$ we get that
\begin{equation*}
T_{\ep}\lesssim\min\big\{\ep^{-F^{-1}_{SG,1}(n,p,q)}, \ep^{-F^{-1}_{SG,2}(n,p,q)}\big\}\ .
\end{equation*}
So we come to the estimate
\begin{align*}
T_{\ep}\lesssim\ep^{-\Ga^{-1}_{SG}(n,p,q)}.
\end{align*}
We obtain the first lifespan estimation in Theorem \ref{liu3}.
\subsubsection{{\bf Blow up region\ $\Ga_{SG}(n,p,q)=0$}}\ \\

{\bf {(\uppercase\expandafter{\romannumeral3})Test function 3: $\Psi(t,x)=-\pa_{t}(\eta^{2p'}_{M}(t)b_{a}(t,x))$}}

Thirdly, we choose test function $ \Psi=-\pa_{t}(\eta^{2p'}_{M}b_{a})$, where $ a=\frac{n+1}{2}-\frac{1}{p}$. Noting that
 $$\Psi(t,x)=-\pa_{t}(\eta^{2p'}_{M}b_{a})=-\pa_{t}\eta^{2p'}_{M}b_{a}-\eta^{2p'}_{M}\pa_{t}b_{a}\geq -\pa_{t}\eta^{2p'}_{M}b_{a} > 0\ .$$
By substituting $\Psi(t,x)$ in \eqref{l1}, we have
\begin{align*}
&\ep C(u_{0}, u_{1})+M^{-1}\int_{0}^{M}\int_{\R^{n}}|v|^{q}\eta_{M}^{2p^{\prime}}b_{a} dxdt\\  \nonumber
\lesssim&\int^{M}_{0}\int_{\R^{n}}u_{t}\Big((\pa^{2}_{t}\eta^{2p'}_{M})b_{a}+2(\pa_{t}\eta^{2p'}_{M})(\pa_{t}b_{a})-D_{1}(x)(\pa_{t}\eta^{2p'}_{M})b_{a}\Big)dxdt\\
= &I_{27}+I_{28}+I_{29}\nonumber
\end{align*}
We apply H\"older inequality and Lemma \ref{le3} to get that 
\begin{align*} \nonumber
|I_{27}|\lesssim &M^{-2}(\int^{M}_{0}\int_{\R^{n}}|u_{t}|^{p}\theta^{2p'}_{M}dxdt)^{\frac{1}{p}}\Big(\int^{M}_{\frac{M}{2}}\int^{t+R}_{0}(t+R)^{-\frac{n-1}{2}p'}(t+R+1-r)^{-1}r^{n-1}drdt\Big)^{\frac{1}{p'}}\\  
\lesssim &M^{-2}\Big(\int^{M}_{0}\int_{\R^{n}}|u_{t}|^{p}\theta^{2p'}_{M}dxdt\Big)^{\frac{1}{p}}(\int^{M}_{\frac{M}{2}}\int^{t+R}_{0}(t+R+1)^{\frac{(n-1)(p-2)}{2p-2}}(t+R+1-r)^{-1}drdt)^{\frac{1}{p'}}\\ 
\lesssim &M^{\frac{np-2n-3p}{2p}}(\ln M)^{\frac{p-1}{p}}\Big(\int^{M}_{0}\int_{\R^{n}}|u_{t}|^{p}\theta^{2p'}_{M}dxdt\Big)^{\frac{1}{p}} ,\nonumber
\end{align*}
\begin{align*} \nonumber
|I_{28}|\lesssim &M^{-1}\Big(\int^{M}_{0}\int_{\R^{n}}|u_{t}|^{p}\theta^{2p'}_{M}dxdt\Big)^{\frac{1}{p}}\Big(\int^{M}_{\frac{M}{2}}\int^{t+R}_{0}(t+R)^{-\frac{n-1}{2}p'}r^{n-1}(t+R+1-r)^{\frac{1}{p}-2}drdt\Big)^{\frac{1}{p'}}\\ \nonumber
\lesssim &M^{-1}\Big(\int^{M}_{\frac{M}{2}}(t+R)^{n-1-\frac{n-1}{2}p'}\ln(t+R)\Big)^{\frac{1}{p'}}\Big(\int^{M}_{0}\int_{\R^{n}}|u_{t}|^{p}\theta^{2p'}_{M}dxdt\Big)^{\frac{1}{p}}\\ 
\lesssim &M^{\frac{np-2n-3p}{2p}}(\ln M)^{\frac{p-1}{p}}\Big(\int^{M}_{0}\int_{\R^{n}}|u_{t}|^{p}\theta^{2p'}_{M}dxdt\Big)^{\frac{1}{p}},\ 
 \end{align*}

\begin{align*}  \nonumber
|I_{29}|\lesssim &M^{-1}\Big(\int^{M}_{0}\int_{\R^{n}}|u_{t}|^{p}\theta^{2p'}_{M}dxdt\Big)^{\frac{1}{p}}\Big(\int^{M}_{\frac{M}{2}}\int^{t+R}_{0}(t+R)^{-\frac{n-1}{2}p'}r^{n-1}(t+R+1-r)^{\frac{1}{p}-2}(1+r)^{-\be p'}drdt\Big)^{\frac{1}{p'}}\\  \nonumber
\lesssim &M^{-1}\Big(\int^{M}_{\frac{M}{2}}(t+R+1)^{-\frac{n-1}{2}p'}(t+R+1-r)^{-1}r^{n-1}(1+r)^{-\be p'}drdt\Big)^{\frac{1}{p'}}\Big(\int^{M}_{0}\int_{\R^{n}}|u_{t}|^{p}\theta^{2p'}_{M}dxdt\Big)^{\frac{1}{p}}\\ 
\lesssim &M^{\frac{np-2n-p}{2p}-\be}\Big(\int^{M}_{0}\int_{\R^{n}}|u_{t}|^{p}\theta^{2p'}_{M}dxdt\Big)^{\frac{1}{p}}\\ 
\lesssim&M^{\frac{np-2n-3p}{2p}}(\ln M)^{\frac{p-1}{p}}\Big(\int^{M}_{0}\int_{\R^{n}}|u_{t}|^{p}\theta^{2p'}_{M}dxdt\Big)^{\frac{1}{p}}. \nonumber
\end{align*}
Therefor we conclude from the estimates $I_{27}$-$I_{29}$
\begin{align}
\label{l20}
\int_{0}^{M}\int_{\R^{n}}|v|^{q}\eta_{M}^{2p^{\prime}}b_{a} dxdt\lesssim M^{\frac{np-2n-p}{2p}}(\ln M)^{\frac{p-1}{p}}\Big(\int^{M}_{0}\int_{\R^{n}}|u_{t}|^{p}\theta^{2p'}_{M}dxdt\Big)^{\frac{1}{p}} .
\end{align}
For the case $F_{SG,1}(n,p,q)=0>F_{SG,2}(n.p,q).$ Observing that the condition $F_{SG,1}(n,p,q)=0$ yields
\begin{equation*}
\Big(n-\frac{n-1}{2}p-(\frac{n-1}{2}-\frac{1}{q})\Big)q+\Big(n-\frac{n-1}{2}q-(\frac{n-1}{2}-\frac{1}{p})-1\Big)=-F_{SG,1}(n,p,q)=0.
\end{equation*}
Taking \eqref{l17} into \eqref{l20} with $ F_{SG,1}(n,p,q)=0$ and by the lower bound of $b_{a}$ in Lemma \ref{le3}, we have that
\begin{align}
\label{l21}
\Big(\int_{0}^{M}\int_{\R^{n}}|v|^{q}\eta_{M}^{2p^{\prime}}b_{a} dxdt\Big)^{pq}\lesssim (\ln M)^{q(p-1)}\Big(\int^{M}_{0}\int_{\R^{n}}|v|^{q}\theta^{2p'}_{M}b_{a}dxdt\Big),
\end{align}
By combining \eqref{l17}, \eqref{l18} and Lemma \ref{le3}, we acquire that
\begin{align}
\label{l22}
\ep^{pq}\lesssim\int_{0}^{M}\int_{\R^{n}}|v|^{q}\theta_{M}^{2p^{\prime}}b_{a} dxdt.
\end{align}
From \eqref{lp1}, \eqref{lp2}, \eqref{l21} and \eqref{l22}, we can get that
\begin{align*}
\ep^{pq}\lesssim& MY'[b_{a}|v|^{q}](M),\\
[MY[b_{a}|v|^{q}](M)]^{pq}\lesssim&(\ln M)^{q(p-1)}MY'[b_{a}|v|^{q}](M).
\end{align*}
Exploiting Lemma \ref{le4} with $p_{2}=q(p-1)+1,$ $p_{1}=pq,$ $\delta=\ep^{pq},$ which lead to
\begin{align*}
T_{\ep}\lesssim \exp(\ep^{-p(pq-1)}).
\end{align*}

{\bf {(\uppercase\expandafter{\romannumeral4})Test function 4: $\Psi(t,x)=\eta^{2p'}_{M}(t)b_{a}(t,x)$}}

We choose test function $\Psi=\eta^{2p'}_{M}b_{a}$, where $a=\frac{n-1}{2}-\frac{1}{q}$. By substituting $\Psi(t,x)$ in \eqref{l2}, we get that
\begin{align} \nonumber
&\ep C(v_{0},v_{1})+\int_{0}^{M}\int_{\R^{n}}|u_{t}|^{p} \eta^{2p'}_{M}b_{a}dxdt\\ 
&=\int^{M}_{0}\int_{\R^{n}}v\Big((\pa^{2}_{t}\eta^{2p'}_{M})b_{a}+2(\pa_{t}\eta^{2p'}_{M})(\pa_{t}b_{a})-D_{2}(x)(\pa_{t}\eta^{2p'}_{M}b_{a})\Big)dxdt\\
&= I_{30}+I_{31}+I_{32}\nonumber
\end{align}
For the case $F_{SG,2}(n,p,q)=0>F_{SG,1}(n.p,q)$. Observing that the condition $F_{SG,2}(n,p,q)=0$ yields
\begin{equation*}
\Big(n-\frac{n-1}{2}p-(\frac{n-1}{2}-\frac{1}{q})\Big)+\Big(n-\frac{n-1}{2}q-(\frac{n-1}{2}-\frac{1}{p})-1\Big)p=-F_{SG,2}(n,p,q)=0.
\end{equation*}
Therefor we conclude from the estimates $I_{30}$-$I_{32}$ that
\begin{align}
\label{l23}
\Big(\int^{M}_{0}\int_{\R^{n}}|u_{t}|^{p}\eta^{2p'}_{M}b_{a}dxdt\Big)^{pq}\lesssim(\ln M)^{p(q-1)}\Big(\int^{M}_{0}\int_{\R^{n}}|u_{t}|^{p}\theta^{2p'}_{M}b_{a}dxdt\Big).
\end{align}
By combining \eqref{l16}, \eqref{l19} and  Lemma \ref{le3}, we acquire
\begin{align}
\label{l24}
\ep^{pq}\lesssim\int^{M}_{0}\int_{\R^{n}}|u_{t}|^{p}\theta^{2p'}_{M}b_{a}dxdt.
\end{align} 
From \eqref{lp1}, \eqref{lp2}, \eqref{l23} and \eqref{l24}, we can get that
\begin{align*}
\ep^{pq}\lesssim& MY'[b_{a}|u_{t}|^{p}](M),\\
[MY[b_{a}|u_{t}|^{p}](M)]^{pq}\lesssim&(\ln M)^{p(q-1)}MY'[b_{a}|u_{t}|^{p}](M).
\end{align*}
Exploiting Lemma \ref{le4} with $p_{2}=p(q-1)+1,$ $p_{1}=pq,$ $\delta=\ep^{pq},$ which lead to
\begin{align*}
T_{\ep}\lesssim \exp(\ep^{-q(pq-1)}).
\end{align*}
For the case $F_{SG,2}(n,p,q)=0=F_{SG,1}(n,p,q)$, that is 
\begin{align*}
n-\frac{n-1}{2}p=\frac{n-1}{2}-\frac{1}{q},\ \ n-\frac{n-1}{2}q=\frac{n+1}{2}-\frac{1}{p}.
\end{align*}
By combining \eqref{l18} and the lower bound of $b_{a}$ in Lemma \ref{le3}, we have that
\begin{align}
\label{lp10}
\ep^{p}\lesssim\int^{M}_{0}\int_{\R^{n}}|u_{t}|^{p}b_{a}\theta^{2p'}_{M}dxdt.
\end{align}
From \eqref{lp1}, \eqref{lp2}, \eqref{l23} and \eqref{lp10}, we can get that
\begin{align*}
\ep^{p}\lesssim& MY'[b_{a}|u_{t}|^{p}](M),\\
[MY[b_{a}|u_{t}|^{p}](M)]^{pq}\lesssim&(\ln M)^{p(q-1)}MY'[b_{a}|u_{t}|^{p}](M).
\end{align*}
Exploiting Lemma \ref{le4} with $p_{2}=p(q-1)+1,$ $p_{1}=pq,$ $\delta=\ep^{p},$ which lead to
\begin{align*}
T_{\ep}\lesssim \exp(\ep^{-(pq-1)}).
\end{align*}
We complete the proof of Theorem \ref{liu3}.

\bibliographystyle{plain1}

\end{document}